\documentclass{amsart}

\theoremstyle{definition}

\usepackage{amsmath}
\usepackage{centernot}
\usepackage{amssymb,latexsym}
\usepackage{enumerate}
\usepackage{cases}
\usepackage{mathtools}

\theoremstyle{remark}

\numberwithin{equation}{section}

\begin{document}

\baselineskip=16pt

\title{Distribution of neighboring values of the Liouville and 
M\"obius functions}

\author{Qi Luo and Yangbo Ye}

\address{Qi Luo: qi-luo@uiowa.edu}

\address{Yangbo Ye: yangbo-ye@uiowa.edu}

\address{Department of Mathematics, The University of 
Iowa, Iowa City, Iowa 52242, USA}



\begin{abstract}
Let $\lambda(n)$ and $\mu(n)$ denote the Liouville function 
and the M\"obius function, respectively. In this study, 
relationships between the values of $\lambda(n)$ and 
$\lambda(n+h)$ up to $n\leq10^8$ for $1\leq h\leq1,000$ are 
explored. Chowla's conjecture predicts that the conditional 
expectation of $\lambda(n+h)$ given $\lambda(n)=1$ for 
$1\leq n\leq X$ converges to the conditional expectation of 
$\lambda(n+h)$ given $\lambda(n)=-1$ for $1\leq n\leq X$ 
as $X\rightarrow\infty$. However, for finite $X$, 
these conditional expectations are different. The observed 
difference, together with the significant difference in 
$\chi^2$ tests of independence, reveals hidden additive 
properties among the values of the Liouville function. 
Similarly, such additive structures for $\mu(n)$ for 
square-free $n$'s are identified. These findings 
pave the way for developing possible, and hopefully efficient, 
additive algorithms for these functions. The potential 
existence of fast, additive algorithms for $\lambda(n)$ and 
$\mu(n)$ may eventually provide scientific evidence 
supporting the belief that prime factorization of large 
integers should not be too difficult. For $1\leq h\leq1,000$, 
the study also tested the convergence speeds of Chowla's 
conjecture and found no relation on $h$. 
\end{abstract}

\maketitle

\section{Introduction} 

Public-key cryptography has become an essential component 
of communication, internet security, digital finance, and 
daily life. Its significance has surpassed the original 
vision set forth in the RSA paper \cite{RivShaAdl}. However, 
the security of public-key cryptography relies primarily 
on the belief that integers cannot be factored into a 
product of primes in polynomial time of the bit length of 
the integers, a belief commonly known as ``not in class P''. 
There is no theoretical or scientific evidence supporting 
this belief, and an increasing number of number theorists 
(cf. Sarnak \cite{Sar3Lect}) are beginning to question 
its validity.

The Liouville function $\lambda(n)=(-1)^{\Omega(n)}$, 
where $\Omega(n)$ is the number of prime factors of $n$ 
counting multiplicities, provides only the parity of the 
number of prime factors of $n$ and, therefore, contains 
much less information than the full prime factorization 
of $n$. Despite this limitation, all known algorithms for 
computing $\lambda(n)$ rely on prime factorization and 
are, therefore, not in class P.

The M\"obius function $\mu(n)$ is defined as 
$\mu(n)=\lambda(n)$ if $n$ is square-free, and $=0$ 
otherwise. Thus, $\mu(n)$ indicates whether $n$ has an 
even or odd number of prime factors when $n$ is 
square-free, while $\mu^2(n)$ indicates whether $n$ is 
square-free or not. Previous algorithms for computing 
$\mu(n)$ and $\mu^2(n)$ have relied on factorization or 
other techniques (such as the work of Booker, Hiary, 
and Keating \cite{BooHiaKea} for $\mu^2(n)$), but they 
have not been shown to be in class P.

The computational complexity and possible algorithms in 
class P of $\lambda(n)$, $\mu(n)$, and $\mu^2(n)$ are 
seemingly easier problems than those of prime 
factorization. These problems have become central 
research topics in number theory and cryptography with 
efforts being mostly focused on the randomness and 
dynamics of these functions. 

In search for possible efficient algorithms, we will take 
another approach. We observe that 
although finding divisors of an integer is not known in P, 
finding the greatest common divisor of two integers is in 
P using the Euclidean algorithm. A notable feature of the 
Euclidean algorithm is that it is an additive algorithm, 
as opposite to the multiplicative algorithms used for 
prime factorization. Most efficient known algorithms of 
factorization use sieve methods which are also additive 
techniques. These suggest that efficient algorithms for 
$\lambda(n)$, $\mu(n)$, and $\mu^2(n)$, 
if they exist, might be additive algorithms. 

With this outlook in mind, the present paper will explore 
additive properties of $\lambda(n)$, $\mu(n)$, and 
$\mu^2(n)$. More precisely, we will get the results 
in the following sections.

\begin{description}

\item[\S2] The values of the squared M\"obius function 
$\mu^2(n)$ have proven additive relationships. 

\item[\S3] The Liouville function $\lambda(n)$.

\begin{description}

\item[\S3.1] Chowla's conjecture for $f(n)=n(n+h)$ 
\cite[(341)]{Cho} predicts 
that for fixed $h\neq0$, the conditional expectation of 
$\lambda(n+h)$ on $\lambda(n)=1$ for $1\leq n\leq X$ 
converges to the conditional expectation of 
$\lambda(n+h)$ on $\lambda(n)=-1$ for $1\leq n\leq X$, 
as $X\rightarrow\infty$. For finite $X$, however, these 
conditional expectations do differ, and hence exhibit 
an additive relationship. For $1\leq h\leq1,000$ and 
$X$ up to $10^8$, this phenomenon will be observed  
numerically. 

\item[\S3.2] The convergent speeds in \S3.1 appear 
to follow the square-root saving numerically. The 
independence of the convergent speeds on $h$ will 
be tested statistically. 

\item[\S3.3] We compare the proportion of 
$\lambda(n+h)=1$ under $\lambda(n)=1$ with the proportion 
of $\lambda(n+h)=1$ under $\lambda(n)=-1$. The $\chi^2$ 
tests will not reject their independence for 
most $1\leq h\leq1,000$. Since values of the Liouville 
function are not random but deterministic, these  
$\chi^2$ tests for various $h$ are themselves indicators 
of relationship between $\lambda(n+h)$ and $\lambda(n)$. 

\end{description}

\item[\S4] If one replaces the Liouville function by the 
M\"obius function, then Chowla's conjecture becomes 
prediction on corresponding sums on 
square-free $n$ with $n+h$ also being square-free. Similar 
computation leads to the same conclusions as in \S\S3.1-3.3.

\end{description}

\noindent
These additive relationships on values of $\lambda(n)$ and 
values of $\mu(n)$ provide a scientific foundation for a 
novel additive algorithm of $\mu(n)$ by Qin and Ye 
\cite{QinYe}.

The techniques used in this paper include the sieve of 
Eratosthenes to generate a database of values of 
$\lambda(n)$ and $\mu(n)$ for $n$ up to $10^8+10^3$. 

\setcounter{equation}{0}
\section{The squared M\"obius function}

According to the work by Carlitz \cite{Car}, Hall \cite{Hal}, 
Heath-Brown \cite{Hea}, Tsang \cite[Theorem 1]{Tsa}, etc, 
we know an asymptotic formula
\begin{equation}\label{SFasymp}
\frac1X
\sum_{n\leq X}
\mu^2(n+h_1)
\cdots
\mu^2(n+h_r)
=
A(h_1,\ldots,h_r)
+
O\big(X^{-\frac13}\big)
\end{equation}
for fixed distinct integers $h_1,\ldots,h_r$, where 
$$
A(h_1,\ldots,h_r)
=
\prod_p
\Big(
1-\frac{u(p)}{p^2}
\Big),
$$
with $u(p)$ being the number of distinct residual classes 
modulo $p^2$ represented by $h_1,\ldots,h_r$. 

The definition of $A(h_1,\ldots,h_r)$ manifests an additive 
relation among $\mu^2(n+h_1),\ldots,\mu^2(n+h_r)$. In 
particular for $r=2$, $h_1=0$, and $h_2=h\geq1$, the 
left hand side of \eqref{SFasymp} equals 
$\frac1X$ times the number of $1\leq n\leq X$ for which 
$n$ and $n+h$ are both square free. Since there are 
approximately $\frac6{\pi^2}X$ square-free numbers in 
$[1,X]$, the left hand side of \eqref{SFasymp} equals 
approximately 
$$
\frac6{\pi^2}
\frac
{\#\{n\in[1,X]|n\,{\rm and}\,n+h\,{\rm are\,square\,free}\}}
{\#\{n\in[1,X]|n\,{\rm is\,square\,free}\}}
$$
which is $\frac6{\pi^2}$ times the conditional probability 
of $n+h$ being square-free on the condition that $n$ is 
square-free for $1\leq n\leq X$. Thus, \eqref{SFasymp} 
means that this conditional probability approaches 
$\frac{\pi^2}6A(0,h)$ as $X\rightarrow\infty$ and hence 
depends on the number of square divisors of $h$. In the 
case of $h$ being square-free, we have 
$A(0,h)=\prod_p(1-2p^{-2})$.  
Therefore, this conditional probability is 
$$
\frac{\pi^2}6
\prod_p 
\Big(1-\frac2{p^2}\Big)
=\prod_p
\frac{p^2}{p^2-1}
\Big(1-\frac2{p^2}\Big)
=\prod_p
\Big(1-\frac1{p^2-1}\Big),
$$ 
which is still different from the unconditional probability 
$$
\frac6{\pi^2}
=\prod_p
\Big(1-\frac1{p^2}\Big)
$$ 
of $n+h$ being square-free, when $X\rightarrow\infty$. 
This additive relationship on the values of $\mu^2(n)$ 
suggests that there might be additive algorithms to 
detect square-free numbers which hopefully are in P. 

\setcounter{equation}{0}
\section{The Liouville function.} 

Chowla \cite[(341)]{Cho} conjectured that for any 
polynomial $f(x)$ of integer coefficients which is 
not of the form $cg^2(x)$ for some 
$c\in\mathbb Z^\times$ and $g\in\mathbb Z[x]$, 
\begin{equation}\label{Chowla}
\sum_{n\leq X}
\lambda\big(f(n)\big)
=o(X)
\end{equation}
with a famous typo $O(X)$ (cf. Sarnak \cite{Sar2}). 
Since $\lambda(n)$ is a completely multiplicative 
function, if $f(n)$ factors, $\lambda\big(f(n)\big)$ 
factors in the same way. 

For $f(n)=n$, \eqref{Chowla} is equivalent to the 
prime number theorem and hence is known:
\begin{equation}\label{f=n}
\frac1X\sum_{n\leq X}
\lambda(n)
=o(1).
\end{equation}
Denote by $L_i$ 
the number of 
$n\leq X$ with $\lambda(n)=(-1)^i$, $i=0,1$. Then 
\eqref{f=n} can be written as 
\begin{equation}\label{1/2}
\lim_{X\rightarrow\infty}\frac{L_0}X
=
\lim_{X\rightarrow\infty}\frac{L_1}X
=\frac12.
\end{equation}

As for the convergence speed in \eqref{f=n}, following 
Soundararajan \cite{Sou} and Balazard and Anne de Roton 
\cite{BalRot}'s results 
$$
\frac1X
\sum_{n\leq X}
\mu(n)
=
O\Big(
\frac1{\sqrt X}
{\rm exp}(c_1(\log X)^\frac12
(\log\log X)^{\frac52+\varepsilon})
\Big)
$$
under the Riemann Hypothesis, 
Humphries \cite{Hum} proved that 
\begin{equation}\label{SumLambda}
\frac1X
\sum_{n\leq X}
\lambda(n)
=
O\Big(
\frac1{\sqrt X}
{\rm exp}(c_2(\log X)^\frac12
(\log\log X)^{\frac52+\varepsilon})
\Big)
\end{equation}
under the Riemann Hypothesis, for certain positive $c_1$ 
and $c_2$. 
On the other hand, Anderson and Stark \cite{AndSta}, 
Borwein, Ferguson, and Mossinghoff \cite{BorFerMos}, and 
Humphries \cite{Hum} proved that there are infinitely many 
integers $X$ such that 
\begin{equation}\label{>sqrt}
\frac1X
\sum_{n\leq X}\lambda(n)
\geq 
\frac{c_3}{\sqrt X}
\end{equation} 
for a constant 
$c_3>0$. Our computation as in Table 1L shows that the 
values of the sum decrease by a factor of $3.264^{-1}$, 
$5.434^{-1}$, $6.295^{-1}$, and $2.170^{-1}$ when 
$X$ increases from $10^4$ to $10^5$ through $10^8$, 
respectively. These rates are around the rate of 
$10^{-\frac12}=3.162^{-1}$ and hence 
confirm the convergence speed as described in 
\eqref{SumLambda} and \eqref{>sqrt}. 

\vspace{4mm}

\begin{center}
\begin{tabular}{|l|c|c|c|c|c|}
\hline
& $X=10^4$ & $X=10^5$ & $X=10^6$ & $X=10^7$ & $X=10^8$
\\
\hline
$\frac1X\sum_{n\leq X}\lambda(n)$
& -0.0094 & -0.00288 & -0.00053 & -8.42E-05 & -3.88E-05
\\
\hline
\end{tabular}

\vspace{2mm}

{\tt Table 1L. Values of $\frac1X\sum_{n\leq X}\lambda(n)$.}

\end{center}

Now let us turn to the case of ${\rm deg}(f)\geq2$ in 
\eqref{Chowla}, which as pointed out by Chowla, ``seems 
an extremely hard conjecture.'' We will 
give numerical evidence for Chowla's conjecture for 
certain polynomials of degrees 2 in \S3.1 below. 

{\it 3.1.} Chowla's conjecture for $f(n)=n(n+h)$.

When $f(n)=n(n+h)$ for a non-zero integer $h$, 
\eqref{Chowla} becomes 
\begin{equation}\label{n(n+h)}
C=
\frac1X
\sum_{n\leq X}
\lambda(n)\lambda(n+h)
=o(1).
\end{equation}
Let $L_{ij}$ be the number of 
$n\in[1,X]$ satisfying 
$\lambda(n)=(-1)^i$ and $\lambda(n+h)=(-1)^j$ for 
$i,j=0,1$. Then Chowla's conjecture in this case is 
\begin{equation}\label{Ln(n+h)}
\lim_{X\rightarrow\infty}
\frac{L_{00}+L_{11}-L_{01}-L_{10}}X=0.
\end{equation}
By \eqref{1/2} and \eqref{Ln(n+h)}, we have
\begin{equation}\label{LimitCondExp}
\frac{L_{00}-L_{01}}{L_0}
-
\frac{L_{10}-L_{11}}{L_1}
=
\frac{L_{00}+L_{11}-L_{01}-L_{10}}X
\cdot
\frac X{L_0}
+
\frac{L_{10}-L_{11}}X
\Big(\frac X{L_0}-\frac X{L_1}\Big)
\rightarrow0.
\end{equation}

Note that 
$$
\frac{L_{00}-L_{01}}{L_0}
=
\frac
{
\sum_{n\leq X,\,\lambda(n)=1}
\lambda(n+h)
}
{
\sum_{n\leq X,\,\lambda(n)=1}1
}
$$
is the conditional expectation of 
$\lambda(n+h)$ on $\lambda(n)=1$ for $1\leq n\leq X$,
while
$$
\frac{L_{10}-L_{11}}{L_1}
=
\frac
{
\sum_{n\leq X,\,\lambda(n)=-1}
\lambda(n+h)
}
{
\sum_{n\leq X,\,\lambda(n)=-1}1
}
$$
is the conditional expectation of 
$\lambda(n+h)$ on $\lambda(n)=-1$ for $1\leq n\leq X$. 
Consequently, Chowla's conjecture as in 
\eqref{LimitCondExp} predicts that these two 
conditional expectations converge to each other 
as $X\rightarrow\infty$. In other words, these 
conditional expectations of $\lambda(n+h)$ approach 
the unconditional expectation. 

Table 2L summarizes our computation which provides a strong 
evidence for Chowla's conjecture, while the full table 
is in the data supplement. 

\vspace{4mm}

\begin{center}
\begin{tabular}{|c|c|c|c|c|c|}
\hline
h & $X=10^4$ & $X=10^5$ & $X=10^6$ & $X=10^7$ & $X=10^8$
\\
\hline
1 & 0.0112 & 0.00068 & -0.00111 & -0.000205 & -3.92E-05\\
2 & 0.0012 & 0.00258 & 6.80E-05 & 0.000125 & 4.63E-05\\
3 & -0.0038 & -0.00074 & -0.000424 & -0.000318 & 0.000107\\
4 & -0.0038 & 0.0013 & -0.000706 & 7.78E-05 & -1.83E-05\\
5 & 0.006 & -0.00176 & 0.000132 & -0.000209 & 1.68E-05\\
10 & -0.0014 & -0.0002 & 0.000102 & -0.000690 & -5.11E-05\\
100 & 0.0022 & 0.00252 & 0.000216 & 0.000152 & -2.95E-05\\
1000 & -0.0098 & -0.00412 & -0.00128 & 9.98E-05 & 0.000121\\
\hline
\end{tabular}

\vspace{2mm}

{\tt Table 2L. Values of 
$C=\frac1X\sum_{n\leq X}\lambda(n)\lambda(n+h)$.}

\end{center}

Consequently, Chowla's conjecture \eqref{n(n+h)} shows 
that the expectation of $\lambda(n+h)$ is independent 
of the value of $\lambda(n)$ but only when 
$X\rightarrow\infty$. For the sake of computational 
complexity and efficient algorithms, the $X$ is finite 
and hence the non-zero entries in Table 2L demonstrate a 
dependence the expectation of $\lambda(n+h)$ on 
$\lambda(n)$. This dependence, however subtle for 
large $X$, represents an additive relation among 
the values of $\lambda(n)$ as what we are looking for. 

{\it 3.2.} Convergence speeds.

Matom\"aki, Radziwi\l\l, and Tao \cite{MatRadTao2015} 
proved an averaged form of \eqref{n(n+h)}:
\begin{equation}\label{Aven(n+h)}
\sum_{h\leq H}
\Big|
\sum_{n\leq X}
\lambda(n)\lambda(n+h)
\Big|
=o(HX)
\end{equation}
as $X\rightarrow\infty$, whenever $H=H(X)\leq X$ goes to 
$\infty$ together with $X$. Their result in \eqref{Aven(n+h)} 
raises a 
question on the convergent speed of \eqref{n(n+h)} on $h$. 
Our numerical computations in this direction are summarized 
in Table 3L. 

\vspace{4mm}

\begin{center}
\begin{tabular}{|l|c|c|c|c|c|}
\hline
 & $X=10^4$ & $X=10^5$ & $X=10^6$ & $X=10^7$ & $X=10^8$
\\
\hline
(a) Absolute values
&&&&&
\\
Mean
& 0.00746 & 0.00250 & 0.000757 & 0.000240 & 7.96E-05
\\
Maximum
& 0.0414 & 0.0113 & 0.00334 & 0.00100 & 0.000324
\\
\hline
(b) Linear regression against $h$
&&&&&
\\
Intercept $b$
& 0.000685 & 0.000225 & -1.29E-05 & -3.50E-05 & -2.22E-05
\\
Slope $m$
& -1.51E-06 & -2.45E-07 & 1.64E-08 & 5.25E-08 & 2.67E-08
\\
R square
& 0.00202 & 0.000514 & 2.47E-05 & 0.00254 & 0.00619
\\
\hline
(c) Correlation with $h$
& -0.0450 & -0.0227 & 0.00497 & 0.0504 & 0.0787
\\
\hline
\end{tabular}

\vspace{2mm}

{\tt Table 3L. Statistics of the values of 
$C=\frac1X\sum_{n\leq X}\lambda(n)\lambda(n+h)$ 
for $1\leq h\leq1,000$: 
(a) Mean and maximum of $|C|$. 
(b) Linear model $C=mh+b+\varepsilon$.
(c) Correlation of $C$ and $h$.}

\end{center}

Part (a) of Table 3L tests the mean and maximum of the 
absolute values of $C$ in \eqref{n(n+h)} over 
$1\leq h\leq1,000$ for $X=10^k$, $k=4,\ldots,8$. The 
pattern of the results fits the right hand side of 
\eqref{SumLambda} and \eqref{>sqrt} and Table 1L well 
and hence suggests that \eqref{SumLambda} and 
\eqref{>sqrt} might hold for \eqref{n(n+h)} as well. 

Part (b) of Table 3L tests a possible linear dependence of 
$C$ on $h$. From the small values of Slope $m$ and R square, 
dependence on $h$ can be rejected. The small correlations 
of $C$ and $h$ in Part (c) provide further support to this 
conclusion. 

{\it 3.3.} The $\chi^2$ test of independence.

Recall that $L_{ij}$ is the number of 
$n\in[1,X]$ satisfying 
$\lambda(n)=(-1)^i$ and $\lambda(n+h)=(-1)^j$ for 
$i,j=0,1$. Denote 
$L_{i+}=L_{i0}+L_{i1}$,  
$L_{+j}=L_{0j}+L_{1j}$, 
$\ell_{ij}=L_{ij}/X$, 
$\ell_{i+}=L_{i+}/X$, and 
$\ell_{+j}=L_{+j}/X$. Then 
$$
\sum_{i=0}^1
\sum_{j=0}^1
L_{ij}
=X
,\ \ \ 
\sum_{i=0}^1
\sum_{j=0}^1
\ell_{ij}
=1.
$$
The independence of $\lambda(n)$ on $\lambda(n+h)$ can be 
formulated as the null hypothesis 
$$
H_0:\ \ 
\ell_{ij}=\ell_{i+}\ell_{+j} 
\ \ {\rm for}\ i,j=0,1. 
$$
The alternative 
hypothesis is $H_1$: the hypothesis $H_0$ is not true; 
that is, $\lambda(n)$ on $\lambda(n+h)$ are not independent. 
The $\chi^2$ test of independence in the contingency table 
$\big(L_{ij}\big)_{i,j=0,1}$ uses the $\chi^2$ statistic 
$$
Q=
\sum_{i=0}^1
\sum_{j=0}^1
\frac{(L_{ij}-\hat E_{ij})^2}{\hat E_{ij}}
$$
with 1 degree of freedom, 
where $\hat E_{ij}=L_{i+}L_{+j}/X$ is the maximum 
likelihood estimate of expected $L_{ij}$ 
(cf. DeGroot and Schervish \cite[\S9.3]{DeGSch}). 
Note that $Q$ also tests the 
homogeneity of the conditional probabilities $\ell_{ij}$ 
(cf. \cite[\S9.4]{DeGSch}). 
Our computation on the $\chi^2$ test is summarized in 
Table 4L with a few entries of large $Q$'s. The full 
table is in the data supplement. 

\vspace{4mm}

\begin{center}
\begin{tabular}{|l|c|c|c|c|c|}
\hline
h & $X=10^4$ & $X=10^5$ & $X=10^6$ & $X=10^7$ & $X=10^8$
\\
\hline
1 & 1.23490 & 0.04512 & 1.22829 & 0.41946 & 0.15336\\
2 & 0.01236 & 0.66141 & 0.00459 & 0.15523 & 0.21454\\
3 & 0.15107 & 0.05599 & 0.18001 & 1.01128 & 1.14787\\
4 & 0.15107 & 0.16687 & 0.49883 & 0.06052 & 0.03364\\
5 & 0.34998 & 0.31265 & 0.01735 & 0.43684 & 0.02808\\
10 & 0.02210 & 0.00433 & 0.01035 & 4.75834 & 0.26134\\
100 & 0.04484 & 0.63111 & 0.04654 & 0.23224 & 0.08692\\
107 & 0.12131 & 1.23363 & 0.01857 & 6.95901 & 4.48429\\
391 & 0.72065 & 3.03895 & 3.59377 & 2.20894 & 3.85567\\
760 & 0.01148 & 4.88087 & 0.12269 & 3.68700 & 9.38083\\
923 & 1.10632 & 0.54436 & 0.03470 & 4.80258 & 4.64668\\
1000 & 0.97609 & 1.70345 & 1.62886 & 0.09959 & 1.45488\\
\hline
\end{tabular}

\vspace{2mm}

{\tt Table 4L. Values of the $\chi^2$ test $Q$ of the 
contingency table for $\lambda(n)$ and $\lambda(n+h)$.}

\end{center}

Hypothesis $H_0$ can be rejected with a $95\%$ confidence 
level if $Q>3.84146$. Our computation in Table 4L shows 
that independence of $\lambda(n)$ and $\lambda(n+h)$ cannot 
be rejected for almost all $h$ and $X$. 

Values of $\lambda(n)$ and $\lambda(n+h)$ are, however, not 
random variables. They form well-defined, deterministic 
sequences. For random variables, the $\chi^2$ test scores 
$Q$ don't have much meaning other than rejecting 
Hypothesis $H_0$ when $Q>$ a certain $c$ under a certain 
confidence level, because one will get a different $Q$ 
using a different sample of the random variable. For 
deterministic sequences, $Q$ scores don't change for 
a given range and hence are intrinsic to the nature 
of the sequences. In particular in our case, the $Q$ 
scores in Table 4L represent properties between 
$\lambda(n)$ and $\lambda(n+h)$ for various $h$ and $X$. 
The fact that these $Q$ scores are vastly different 
shows that the relationship between $\lambda(n)$ and 
$\lambda(n+h)$ changes with $h$ and $X$. This in turn 
demonstrates an additive relationship among values of 
the Liouville function. 

\setcounter{equation}{0}
\section{The M\"obius function on square-free numbers.} 

The sums $\sum_{n\leq X}\mu(n)$ and 
$\sum_{n\leq X}\mu(n)\mu(n+h)$ are similar to the sums 
in \eqref{f=n} and \eqref{n(n+h)} but with $n$ and 
$n+h$ being square-free. Tables 1M through 4M below 
summarize our computation for the M\"obius function 
which leads to the same conclusions as in \S\S3.1-3.3. 

\vspace{3mm}

\begin{center}
\begin{tabular}{|l|c|c|c|c|c|}
\hline
& $X=10^4$ & $X=10^5$ & $X=10^6$ & $X=10^7$ & $X=10^8$
\\
\hline
$\frac1{Y_1}\sum_{n\leq X}\mu(n)$
& -0.00378 & -0.000790 & 0.000349 & 0.000171 & 3.17E-05
\\
\hline
\end{tabular}

\vspace{2mm}

{\tt Table 1M. Values of $\frac1{Y_1}\sum_{n\leq X}\mu(n)$, 
where $Y_1$ is the number of square-free $n\in[1,X]$.}

\vspace{4mm}

\begin{tabular}{|c|c|c|c|c|c|}
\hline
h & $X=10^4$ & $X=10^5$ & $X=10^6$ & $X=10^7$ & $X=10^8$
\\
\hline
1 & 0.00372 & -0.00580 & 0.00127 & 0.000522 & -8.15E-05\\
2 & -0.00526 & 0.00294 & -0.00119 & 3.38E-05 & -3.07E-05\\
3 & -0.00372 & -0.00353 & -0.000942 & 4.25E-05 & 3.87E-05\\
4 & 0.00371 & -0.000558 & -0.00108 & 0.000189 & -0.000118\\
5 & -0.000309 & -0.00316 & 0.000496 & -0.000424 & -8.92E-05\\
10 & -0.0220 & -0.0104 & -0.00227 & -0.00159 & -0.000242\\
100 & -0.00752 & -0.000436 & -0.000452 & -0.000581 & -4.42E-05\\
1000 & -0.0180 & -0.00596 & -0.00111 & -0.000289 & -4.15E-05\\
\hline
\end{tabular}

\vspace{2mm}

{\tt Table 2M. Values of 
$D=\frac1{Y_2}\sum_{n\leq X}\mu(n)\mu(n+h)$, where $Y_2$ 
is the number of square-free $n\in[1,X]$ such that $n+h$ 
is also square-free.}

\vspace{4mm}

\begin{tabular}{|l|c|c|c|c|c|}
\hline
 & $X=10^4$ & $X=10^5$ & $X=10^6$ & $X=10^7$ & $X=10^8$
\\
\hline
(a) Absolute values
&&&&&
\\
Mean
& 0.0129 & 0.00408 & 0.00124 & 0.000421 & 0.000135
\\
Maximum
& 0.0850 & 0.0185 & 0.00625 & 0.00168 & 0.000555
\\
\hline
(b) Linear regression against $h$
&&&&&
\\
Intercept $b$
& -0.000207 & 0.000341 & -8.64E-05 & -5.88E-05 & -2.13E-05
\\
Slope $m$
& -6.39E-07 & -4.57E-07 & 2.03E-07 & 8.71E-08 & 2.59E-08
\\
R square
& 0.000115 & 0.000655 & 0.00137 & 0.00228 & 0.00193
\\
\hline
(c) Correlation with $h$
& -0.0107 & -0.0256 & 0.0370 & 0.0478 & 0.0439
\\
\hline
\end{tabular}

\vspace{2mm}

{\tt Table 3M. Statistics of the values of 
$D=\frac1{Y_2}\sum_{n\leq X}\mu(n)\mu(n+h)$ 
for $1\leq h\leq1,000$: 
(a) Mean and maximum of $|D|$. 
(b) Linear model $D=mh+b+\varepsilon$.
(c) Correlation of $D$ and $h$.}

\vspace{4mm}

\begin{tabular}{|l|c|c|c|c|c|}
\hline
h & $X=10^4$ & $X=10^5$ & $X=10^6$ & $X=10^7$ & $X=10^8$
\\
\hline
1 & 0.04702 & 1.08525 & 0.51843 & 0.87790 & 0.21455\\
2 & 0.09417 & 0.27883 & 0.45534 & 0.00369 & 0.03050\\
3 & 0.04458 & 0.40277 & 0.28575 & 0.00580 & 0.04834\\
4 & 0.06671 & 0.01509 & 0.56955 & 0.17298 & 0.67732\\
5 & 0.00032 & 0.32247 & 0.07947 & 0.58047 & 0.25687\\
10 & 1.56487 & 3.47856 & 1.65657 & 8.11887 & 1.89446\\
100 & 0.28372 & 0.00955 & 0.10290 & 1.70498 & 0.09883\\
109 & 0.08963 & 0.00920 & 3.16200 & 6.10480 & 9.35062\\
298 & 0.00045 & 0.21368 & 0.37504 & 4.28434 & 7.50830\\
374 & 0.07947 & 1.56817 & 1.24637 & 7.02544 & 4.29383\\
391 & 0.12518 & 1.08537 & 3.34727 & 4.02410 & 8.32421\\
923 & 0.01960 & 0.43887 & 1.08631 & 5.14136 & 5.38660\\
1000 & 1.64232 & 1.79404 & 0.62100 & 0.42049 & 0.08692\\
\hline
\end{tabular}

\vspace{2mm}

{\tt Table 4M. Values of the $\chi^2$ test $Q$ of the 
contingency table for $\mu(n)$ and $\mu(n+h)$ when 
$n$ and $n+h$ are both square-free.}

\end{center}

\setcounter{equation}{0}
\section{Computation procedures.}

The hardware used in this study is a MacBook Pro with an 
M2 chip and 16GB memory. The programing language used 
is Python. A database of $\Omega(n)$, $\lambda(n)$, and 
$\mu^2(n)$ for $n\leq10^8+10^3$ was 
generated by the sieve of Eratosthenes in 500sec. 
The size of this database is 3.2GB as a NumPy array file. 

For $X=10^8$, each entry in Tables 2L, 4L, 2M, and 4M 
requires a run-time of 40sec to 60sec. Since we computed 
these entries for 1,000 $h$'s, the run-times for the 
$X=10^8$ columns in Tables 2L, 4L, 2M, and 4M are 
55,387sec, 61,374sec, 40,906sec, and 50,130sec, 
respectively. 

The reason that we computed 1,000 $h$'s is to have a 
large sample size for statistical analysis on convergent 
speeds in Tables 3L and 3M. If one reduces the number of 
$h$'s, the whole project can be extended to $10^9$ 
within the same run-time frame using the same hardware 
setup. The computation can be further extended using a 
computer cluster as the algorithms are readily parallel. 

Statistical analyses in Tables 3L and 3M were performed 
using Excel. 

\vspace{5mm}

{\bf Acknowledgments.} The first author was supported in 
part by an undergraduate research assistantship in the 
Department of Mathematics, The University of Iowa.

\end{document}